\documentclass[preprint,12pt]{elsarticle}

\usepackage{amssymb}
\usepackage{amsmath}
\biboptions{numbers,sort&compress}

\newtheorem{defn}{Definition}[section]

\renewcommand{\raggedright}{\leftskip=0pt \rightskip=0pt plus 0cm}
\raggedright
\def\bb{\begin}
\def\bc{\begin{center}}
\def\ec{\end{center}}
\def\be{\begin{equation}}
\def\ee{\end{equation}}
\def\ba{\begin{array}}
\def\ea{\end{array}}
\def\bea{\begin{eqnarray}}
\def\eea{\end{eqnarray}}
\def\beaa{\begin{eqnarray*}}
\def\eeaa{\end{eqnarray*}}

\def\lto{\rightarrow}

\def\to{\lto}

\def\ifl{\iffalse}
\def\Proof{\noindent{\bf Proof} \quad}
\def\qed{\hfill $\Box$ \smallskip}

\def\nn{\nonumber}
\def\lb{\label}

\begin{document}

\begin{frontmatter}



\title{Quasi-selfadjoint extensions of dual pairs of linear relations}


%\ifl

\author{
Guixin Xu$^{a}$
, \qquad Guojing Ren$^{b}\footnote{The corresponding author. \\
\indent \; Email addresses: {\tt guixinxu$_-$ds@163.com} (G. Xu),\;{\tt gjren@sdufe.edu.cn} (G. Ren).}$
%E-mail: {\tt zhangmr@tsinghua.edu.cn}.
}

\address{$^{a}$School of Mathematics and Statistics, Beijing Technology and Business University, Beijing, P. R. China;\\
$^{b}$School of Mathematics and Quantitative Economics, Shandong University of Finance and Economics, Jinan, Shandong, P. R. China

}

%\fi

\begin{abstract}
This paper investigates quasi-selfadjoint extensions of dual pairs of linear relations in Hilbert spaces. Some properties of dual pairs of linear relations are given and an Hermitian linear relation associated with a dual pair of linear relations is introduced.
Necessary and sufficient conditions for quasi-selfadjoint extensions of dual pairs of linear relations are obtained.
These results extends the corresponding results for dual pairs of linear operators to dual pairs of linear relations.
\end{abstract}
\smallskip


\end{frontmatter}

%%%\tableofcontents %%%\newpage

\section{Introduction}
\setcounter{equation}{0} \lb{first}

In the classical operator theory,  the operators are discussed in single-valued case (cf., \cite{Kato84,Weidmann80}). Since the adjoint of a non-densely defined operator is multi-valued,  we have always required that the operators are densely defined when  considering their adjoints in the classical operator theory. In 1950, von Neumann introduced linear relations in order to study the adjoint of a non-densely defined linear differential operator \cite{Neumann50}. Since then, more and more multi-valued operators have been founded and they have attracted a lot of attention from mathematicians. For example,  the minimal and maximal operators corresponding to a linear discrete Hamiltonian system or a linear symmetric difference equation are multi-valued and non-densely defined (cf., \cite{Ren14, SS11}), and therefore the classical  theory of linear operator is not available in this cases. So it is necessary and urgent to study some topics about multi-valued linear operators,which are also called linear relations or linear subspaces (cf., \cite{Arens61,Coddington73,Coddington76,DS74}). Throughout the present paper, a linear relation is called an operator if it is single-valued.

Dual pairs of linear operators and their extensions have been studied by many authors due to their widely applications (cf., \cite{BrownG09,BrownH09,BrownK19,BrownM08,Gorbachuk91,Malamud01,Malamud97}). Recall that a pair of operators $A$ and $B$ forms a dual pair $\{A,B\}$ if
\be
\langle A(f),g\rangle=\langle f,B(g)\rangle,\;f\in D(A),\;g\in D(B).\nn
\ee
In \cite{BrownK19}, the authors gave some necessary and sufficient conditions for quasi-selfadjoint extensions of the dual pair $\{A,B\}$ by
constructing the symmetric operator $S$, which is given by
\be
S=\left(
  \begin{array}{cc}
    0 & A \\
    B & 0 \\
  \end{array}
\right).\nn
\ee
In the present paper, inspired by the methods used in \cite{BrownK19}, we give some analogous quasi-selfadjoint extension results about dual pairs of linear relations. Some results obtained in the present paper extend the corresponding results for dual pairs of operators to  dual pairs of  relations (see Remarks \ref{r1} and \ref{r2}).

The rest of the present paper is organized as follows. In Section 2, some notations, basic concepts
and fundamental results about linear relations are introduced. The concept of the dual pairs of linear relations is recalled and  its some fundamental properties are given. In particular, an Hermitian linear relation associated with a dual pair of linear relations is introduced and the relationship between its deficiency  and those of the dual pair is given. In Section 3,  the quasi-selfadjoint extensions of dual pairs of linear relations are studied. In particular,  some sufficient and necessary conditions  for quasi-selfadjoint extensions of dual pairs of linear relations are obtained.

\section{Preliminaries}
\setcounter{equation}{0} \lb{second}
In this section, we shall first recall some basic concepts and give some fundamental results about linear relations, including the concepts of Hermitian and  self-adjoint relations.  Then we shall give several  properties of dual pairs of linear  relations. Finally, we introduce an Hermitian linear relation associated with a dual pair of relations and study its properties.

This section is divided into three subsections.

\subsection{Some notations and basic concepts about linear relations}

Let $X$ be a Hilbert space over the complex  field ${\mathbb C}$ with the inner product $\langle\cdot,\cdot\rangle$.   The product space
 $X^2$ with the following induced inner product, still denoted by $\langle\cdot,\cdot\rangle$ without any
confusion:
\be
\langle(x_1,y_1),(x_2,y_2)\rangle=\langle x_1, x_2\rangle
+\langle y_1,y_2\rangle,\;\;(x_1,y_1),\;(x_2,y_2)\in X^2.\nn
\ee

Any linear subspace $T\subset X^2$ is called a \emph{linear relation (briefly, relation or subspace)} of $X^2$. By $LR(X)$ denote the set of all linear relations of $X^2$.

Let $T\in LR(X)$. The domain $D(T)$, range $R(T)$, and null space $N(T)$ of $T$ are respectively defined by
$$\begin{array}{rrll}
D(T):&=&\{x\in X:\, (x,y)\in T \;{\rm for\; some}\;y\in X\},\\[0.5ex]
R(T):&=&\{y\in X:\, (x,y)\in T \;{\rm for\; some}\;x\in X\},\\[0.5ex]
N(T):&=&\{x\in X:\,(x,0)\in T\}.\\[0.4ex]
\end{array}$$
Further, denote
$$\begin{array}{rrll}
T(x)&:=&\{y\in X:\,(x,y)\in T\},\\[0.5ex]
T^{-1}&:=&\{(y,x)\in X^2:\,(x,y)\in T\}.\\[0.4ex]
\end{array}$$

Let $T\in LR(X)$.  $T$ is said to be {\it closed} if $T=\overline{T}$, where $\overline{T}$ is the closure of $T$ in $X^2$. By $CR(X)$ denoted the set of all the closed linear relations of $X^2$.

Let $T,\,A\in LR(X)$ and $\alpha \in {\mathbb C}$. Define
$$\begin{array}{lll}
\alpha\, T:= \{(x, \alpha\, y)\in X^2 : (x, y) \in T\},\\[0.5ex]
T + A := \{(x, y + z)\in X^2 : (x, y) \in T,\; (x, z)\in A\},\\[0.5ex]
T-\alpha:=\{(x, y - \alpha x)\in X^2 : (x, y) \in T\}.

\end{array}$$
If $T\cap A=\{(0,0)\}$, denote
\be
T\dot{+} A:=\{(x_{1}+x_{2}, y_{1}+y_{2})\in X^2:\, (x_{1},y_{1})\in T, \, (x_{2},y_{2})\in A\}.\nn
\ee
Further,
if $T$ and $A$ are orthogonal, that is,
$\langle(x, y),(u,v)\rangle = 0 $ for all
$(x, y) \in T$ and $(u, v) \in A$, then denote
\be
T\oplus A:=T\dot{+} A.\nn
\ee
The \emph{product} of linear relations $T$ and $A$ is defined by (see \cite{Arens61})
\be
AT:=\{(x,z)\in X^2:\, (x,y)\in T,\;(y,z)\in A\;{\rm for\; some}\;y\in X\}.\nn
\ee
Note that if $A$ and $T$ are both operators, then $AT$ is also an operator.

\begin{defn}\lb{ad} {\rm Let $T\in LR(X)$. The {\it adjoint} of $T$ is defined  by
$$T^{*}:=\{(f,g)\in X^2:\;\langle g,x\rangle=\langle f,y\rangle\;{\rm for\;all}\;(x,y)\in T\}.$$
$T$ is said to be  {\it Hermitian} if $T\subset T^*$, and said to be  {\it self-adjoint}  if $T=T^*$.}
\end{defn}

Let $T\in LR(X)$ and $\lambda\in{\mathbb C}$. Denote
\be
M_{\lambda}(T)=N(T^*-\lambda).\nn
\ee
Then $M_{\lambda}(T)$ is a closed subspace since $T^*$ is closed. If $T$ is Hermitian, then $\dim M_{\lambda}(T)$ is constant in the upper
and lower half-planes by Corollary 2.2 in \cite{Shi12}. Denote $n_{\pm}(T):=\dim M_{\pm i}(T)$ if $T$ is Hermitian.

Let $T\in LR(X)$ be closed. Arens introduced the following important decomposition \cite{Arens61}:
\be
T=T_{s}\oplus T_{\infty}\nn
\ee
where
\be
T_\infty :=\{(0,y)\in X^2: (0,y)\in T\},\;\;
T_s :=T\ominus T_{\infty}.\nn
\ee
It can be easily to verify that $T_{s}$ is an operator, and $T$ is an operator if and only if $T=T_{s}$.
$T_{s}$ and $T_{\infty}$ are called the operator and pure multi-valued parts of $T$, respectively.
In addition, they satisfy the following properties.
\be
D(T_s)=D(T),\;R(T_s)\subset T(0)^\perp,\;T_\infty =\{0\}\times T(0).\nn
\ee

%\bb{lem}\lb{l2.5}{\rm(\cite[p. 26]{DS85})}  Let $X$  be a Hilbert space and $T\in LR(X)$ be self-adjoint. Then $T_s$ and $T_{\infty}$ are self-adjoint relations in $(T(0)^{\perp})^2$ and $T(0)^2$, respectively.
%\end{lem}

\subsection{Dual pairs of linear relations}

In this subsection, we shall recall the concept of the dual pair of linear relations, which can be found in \cite{Hassi05, Malamud02}, and give its some fundamental properties.

\begin{defn}
{\rm  (1) A pair of closed linear relations $A$ and $B$ in $X^2$ is said to form a dual pair $\{A,B\}$ if
\be
\langle g,h\rangle=\langle f,k\rangle,\;\;(f,g)\in A,\;\;(h,k)\in B,\nn
\ee
or equivalently if $A\subset B^*$ ($\Leftrightarrow B\subset A^*$).

(2) A linear relation $\tilde{A}\in CR(X)$ is called a proper extension of a dual pair $\{A,B\}$ if $A\subset \tilde{A}\subset B^*$. The set of all proper extensions of a dual pair $\{A,B\}$ is denoted by Ext$\{A,B\}$.}
\end{defn}

Let $X_{+A}$ be equipped with the space $D(A^*)$ with the inner product
\be
\langle f,g\rangle_{+A}=\langle f,g\rangle+\langle (A^*)_s(f),(A^*)_s(g)\rangle,\;f,g\in D(A^*).\nn
\ee
Then $X_{+A}$ is a Hilbert space by \cite[Theorem 4.3]{SXR18}. By $\oplus_A$ denote the orthogonal sum in $X_{+A}$.

\bb{lem}\lb{l1}
Let $A,B\in CR(X)$ and $\{A,B\}$ be a dual pair with $(A^*)_s|_{D(B)}=B_s$. Then $D(B)$ is a closed subspace in  $X_{+A}$.
\end{lem}

\Proof
Given any sequence $\{f_n\}_{n=1}^{\infty}\subset D(B)$ with $f_n\to f$ in $X_{+A}$. Then $f_n\to f$ and $(A^*)_s(f_n)\to(A^*)_s(f)$ in $X$ as $n\to\infty$. Note that $B$ is closed and $(A^*)_s|_{D(B)}=B_s$. We have that $B_s$ is closed and
\be
B_s(f_n)=(A^*)_s(f_n)\to(A^*)_s(f),\;as\;n\to\infty.
\ee
Therefore, $f\in D(B_s)=D(B)$, i.e., $D(B)$ is closed  in  $X_{+A}$.
\qed

\bb{lem}\lb{l2}
Let $A,B\in CR(X)$ and $\{A,B\}$ be a dual pair with $(A^*)_s|_{D(B)}=B_s$ and $(B^*)_s|_{D(A)}=A_s$. Then the following decompositions hold
\be\lb{e1}
D(A^*)=D(B)\oplus_A N(1+B^*A^*),
\ee
\be\lb{e2}
D(B^*)=D(A)\oplus_B N(1+A^*B^*).
\ee
In particular, if $A^*=B$ (or equivalently, $A=B^*$), then
\be
N(1+B^*A^*)=\{0\}\;(or\; equivalently,N(1+A^*B^*)=\{0\}).\nn
\ee
\end{lem}

\Proof
We only prove \eqref{e1}, and \eqref{e2} can be proved similarly. It follows from Lemma \ref{l1} that $D(B)$ is a closed subspace in  $X_{+A}$. Obviously, $D(B)\oplus_A N(1+B^*A^*)\subset D(A^*)$. Let $g\in D(B)^{\perp}\subset X_{+A}$. Then for any $f\in D(B)$
\be\lb{e3}
\langle f,g\rangle_{+A}=\langle f,g\rangle+\langle (A^*)_s(f),(A^*)_s(g)\rangle=0.
\ee
Let $h=(A^*)_s(g)$. For every $(f,y)\in B\subset A^*$, $y$ can be decomposed as $y=(A^*)_s(f)+y_{\infty}$, where $y_{\infty}\in  A^*(0)$. Note that $h\in R((A^*)_s)\subset A^*(0)^{\perp}$. Then by \eqref{e3}, we have
\be
\langle f,g\rangle+\langle y,h\rangle=\langle f,g\rangle+\langle (A^*)_s(f),h\rangle+\langle y_{\infty},h\rangle=\langle y_{\infty},h\rangle=0,\nn
\ee
which implies that $(h,-g)\in B^*$. Since $(g,h)\in A^*$, we have $(g,-g)\in B^*A^*$, i.e., $g\in N(1+B^*A^*)$.

On the other hand, let $g\in N(1+B^*A^*)$, i.e., $(g,-g)\in B^*A^*$. There exists $h\in X$ such that $(g,h)\in A^*$ and $(h,-g)\in B^*$. Then $h$ can be decomposed as $h=(A^*)_s(g)+h_{\infty}$ with $h_{\infty}\in A^*(0)$. For every $f\in D(B)\subset X_{+A}$, we have $(A^*)_s(f)=B_s(f)$  by the assumption $(A^*)_s|_{D(B)}=B_s$.  Note that $(A^*)_s(f)\in A(0)^{\perp}$ and $(h,-g)\in B^*$, we have
\beaa
&&\langle f,g\rangle_{+A}=\langle f,g\rangle+\langle (A^*)_s(f),(A^*)_s(g)\rangle\\
&=&\langle f,g\rangle+\langle (A^*)_s(f),h-h_{\infty}\rangle=\langle f,g\rangle+\langle (A^*)_s(f),h\rangle\\
&=&\langle f,g\rangle+\langle B_s(f),h\rangle=0,
\eeaa
which yields that $g\in D(B)^{\perp}$ in $X_{+A}$. Thus, \eqref{e1} holds. This completes the proof.
\qed

Let $A,B\in CR(X)$ and $\{A,B\}$ be a dual pair with $(A^*)_s|_{D(B)}=B_s$ and $(B^*)_s|_{D(A)}=A_s$. Define
\be
n(A^*,B):=\dim(D(A^*)/D(B)),\;n(B^*,A):=\dim(D(B^*)/D(A)).\nn
\ee

Next, we  shall show that the numbers $n(A^*,B)$ and $n(B^*,A)$ are equal under some conditions.

\bb{lem}\lb{l12}
Let $A,B\in CR(X)$ and $\{A,B\}$ be a dual pair with $B^*(0)\cap N(A^*)=\{0\}$ and $A^*(0)\cap N(B^*)=\{0\}$. Then there is a isomorphic mapping between the kernel spaces $N(1+B^*A^*)$ and $N(1+A^*B^*)$.
\end{lem}

\Proof The proof is divided into three steps.

{\bf Step 1.} We show that for every $g\in N(1+B^*A^*)$, there exists only one $h\in N(1+A^*B^*)$ such that
\be\lb{ec1}
(g,h)\in A^*,\;(h,-g)\in B^*.
\ee
For the existence, let $g\in N(1+B^*A^*)$. Then  $(g,-g)\in B^*A^*$, and hence, there is $h\in X$ such that \eqref{ec1} holds, which implies that $(h,-h)\in A^*B^*$, i.e., $h\in N(1+A^*B^*)$.

Now, we prove the uniqueness. If there exists $h_0\in N(1+A^*B^*)$ such that
\be
(g,h_0)\in A^*,\;(h_0,-g)\in B^*,\nn
\ee
then we have that $(0,h-h_0)\in A^*$ and $(h-h_0,0)\in B^*$, and hence, $h-h_0\in A^*(0)\cap N(B^*)=\{0\}$, which implies that $h=h_0$.

Define a mapping $Q:\;N(1+B^*A^*)\to N(1+A^*B^*)$ by
\be\lb{ec2}
(g,Q(g))\in A^*,\;(Q(g),-g)\in B^*, \;g\in N(1+B^*A^*).
\ee

{\bf Step 2.} We show that $Q$ is  bijective. For every $h\in N(1+A^*B^*)$, i.e., $(h,-h)\in A^*B^*$, there exist $g'\in X$ such that
\be\lb{ec3}
(h,g')\in B^*,\;(g',-h)\in A^*.
\ee
Then $(g',-g')\in B^*A^*$, i.e., $g'\in N(1+B^*A^*)$. Let $g=-g'$. Then $g\in N(1+B^*A^*)$ and $Q(g)=h$ by \eqref{ec2} and \eqref{ec3}. This yields that $Q$ is surjective.

On the other hand, if $h=0$ in \eqref{ec3}, then $g'\in B^*(0)\cap N(A^*)=\{0\}$, i.e., $g'=0$. Hence, $g=-g'=0$, which implies that $Q$ is injective. Therefore, $Q$ is bijective.

{\bf Step 3.} We show that $Q$ is  linear. For  any $g_1,g_2\in N(1+B^*A^*)$, we have that
\be
(g_i,Q(g_i))\in A^*,\;(Q(g_i),-g_i)\in B^*,\;i=1,2.\nn
\ee
Then
\beaa
(c_1g_1\pm c_2g_2, c_1Q(g_1)\pm c_2Q(g_2))\in A^*,\\
(c_1Q(g_1)\pm c_2Q(g_2),-(c_1g_1\pm c_2g_2))\in B^*.\nn
\eeaa
Hence, $Q(c_1g_1\pm c_2g_2)=c_1Q(g_1)\pm c_2Q(g_2)$, which implies that $Q$ is linear.

Therefore $Q$ maps $N(1+B^*A^*)$ isomorphically onto $N(1+A^*B^*)$. The proof is complete.
\qed

By Lemma \ref{l12}, one can easily get the following result.

\bb{lem}\lb{l3}
Let $A,B\in CR(X)$ and $\{A,B\}$ be a dual pair with $B^*(0)\cap N(A^*)=\{0\}$ and $A^*(0)\cap N(B^*)=\{0\}$. Then
\be\lb{e6}
\dim N(1+B^*A^*)=\dim N(1+A^*B^*).
\ee
\end{lem}

Lemma \ref{l2} together with \ref{l3} implies the following result hold.

\bb{prop}\lb{p1}
Let $A,B\in CR(X)$ and $\{A,B\}$ be a dual pair with $(A^*)_s|_{D(B)}=B_s$ and $(B^*)_s|_{D(A)}=A_s$. If $B^*(0)\cap N(A^*)=\{0\}$ and $A^*(0)\cap N(B^*)=\{0\}$, then
\be\lb{e5}
n(A^*,B)=n(B^*,A).
\ee
\end{prop}

\subsection{An Hermitian linear relation associated with a dual pair of linear relations}

In this subsection, we shall introduce an Hermitian linear relation associated with a dual pair of linear relations and give its some properties.

Let $A,B\in CR(X)$ and $\{A,B\}$ be a dual pair. Define
\be\lb{S}
S=\left(
    \begin{array}{cc}
      0 & A \\
      B & 0 \\
    \end{array}
  \right)
 \; {\rm with}\; D(S)=D(B)\times D(A),
\ee
i.e., $S=\{(\{x_1,x_2\},\{y_1,y_2\}):\;(x_1,y_2)\in B,\;(x_2,y_1)\in A\}\in LR(X^2)$.

\bb{lem}\lb{l4}
Let $S$ be defined in \eqref{S}. Then
\be\lb{e11}
S^*=\left(
    \begin{array}{cc}
      0 & B^* \\
      A^* & 0 \\
    \end{array}
  \right).
\ee
\end{lem}

\Proof
Note that
\beaa
S^*=\{(\{f_1,f_2\},\{g_1,g_2\}):\langle f_1,y_1\rangle+\langle f_2,y_2\rangle=\langle g_1,x_1\rangle+\langle g_2,x_2\rangle\\
{\rm for\,all}\,(x_1,y_2)\in B,\;(x_2,y_1)\in A\}.
\eeaa
We can easily show that $\left(
    \begin{array}{cc}
      0 & B^* \\
      A^* & 0 \\
    \end{array} \right)\subset S^*$. Let $(\{f_1,f_2\},\{g_1,g_2\})\in S^*$. For any $(x_1,y_2)\in B$, we have that $\langle f_2,y_2\rangle=\langle g_1,x_1\rangle$ by noting that $(0,0)\in A$, which implies that $(f_2,g_1)\in B^*$. Similarly, one can prove that $(f_1,g_2)\in A^*$. Hence, $(\{f_1,f_2\},\{g_1,g_2\})\in \left(
    \begin{array}{cc}
      0 & B^* \\
      A^* & 0 \\
    \end{array} \right)$. Therefore, \eqref{e11} holds.
\qed

\bb{prop}\lb{p2}
Let $A,B\in CR(X)$ and $\{A,B\}$ be a dual pair with $(A^*)_s|_{D(B)}=B_s$ and $(B^*)_s|_{D(A)}=A_s$, and $S$ be defined in \eqref{S}. If $B^*(0)\cap N(A^*)=\{0\}$ and $A^*(0)\cap N(B^*)=\{0\}$,  then
\begin{itemize}
\item[{\rm (i)}] $S$ is an Hermitian relation in $X^2\times X^2$, and
\be\lb{e12}
n_{+}(S)=n_{-}(S)=n(A^*,B)=n(B^*,A).
\ee
\item[{\rm (ii)}] $S$ is self-adjoint if and only if $A=B^*$ (or equivalently $B=A^*$).
\end{itemize}
\end{prop}

\Proof
(i) It follows from Lemma \ref{l4} and the fact that $\{A,B\}$ is a dual pair that
\be
S=\left(
    \begin{array}{cc}
      0 & A \\
      B & 0 \\
    \end{array}
  \right)\subset \left(
                   \begin{array}{cc}
                     0 & B^* \\
                     A^* & 0 \\
                   \end{array}
                 \right)=S^*,\nn
\ee
which implies that $S$ is an Hermitian relation. In order to show \eqref{e12} it suffices to show
\be\lb{e13}
\dim N(S^*\mp i)=\dim N(1+B^*A^*),
\ee
by Lemma \ref{l2}, \eqref{e6} and \eqref{e5}.

Let $f=\{g,h\}\in N(S^*+i)$, i.e., $(f,-if)\in S^*$. By \eqref{e11}, we have that
\be\lb{e14}
(g,-ih)\in A^*,\;(h,-ig)\in B^*,
\ee
which implies that $(g,-g)\in B^*A^*$, i.e., $g\in N(1+B^*A^*)$. Define a projection $P:\;N(S^*+i)\to N(1+B^*A^*)$ by
\be
P(f)=g\;for\;f=\{g,h\}\in N(S^*+i).\nn
\ee
We show that $P$ is a isomorphic mapping from $N(S^*+i)$ onto $N(1+B^*A^*)$.

Given any $g\in N(1+B^*A^*)$, we have that $(g,Q(g))\in A^*$ and $(Q(g),-g)\in B^*$, where $Q$ is defined in \eqref{ec2}. Let $f=\{-ig,Q(g)\}$. Then $(f,-if)\in S^*$, i.e., $f\in N(S^*+i)$. Hence, $P(f)=g$, which yields that $P$ is surjective.

Let $P(f)=0$, where $f=\{0,h\}\in N(S^*+i)$. It follows from \eqref{e14} that $h\in A^*(0)\cap N(B^*)=\{0\}$, i.e., $h=0$, and then $f=0$. This implies that $P$ is injective.

Hence,  the projection $P$ maps $N(S^*+i)$ isomorphically onto $N(1+B^*A^*)$. Consequently, $\dim N(S^*+i)=\dim N(1+B^*A^*)$. Similarly, $\dim N(S^*-i)=\dim N(1+B^*A^*)$.  Therefore, \eqref{e13} holds, which implies that \eqref{e12} holds.

(ii) This statement is derived immediately by (i) and the proof is complete.

\qed

\bb{rem}\lb{r1}
{\rm Proposition \ref{p2} extends \cite[Proposition 2.7]{BrownK19} for dual pairs of operators to dual pairs of linear relations.}
\end{rem}

\section{Quasi-selfadjoint extensions of dual pairs of linear  relations}
\setcounter{equation}{0} \lb{conc}

In this section, we shall introduce the quasi-selfadjoint extensions of dual pairs of linear relations, and give its some sufficient and necessary conditions.

\bb{defn}
{\rm Let $A,B\in CR(X)$ and $\{A,B\}$ be a dual pair.

(1) An extension $\tilde{A}\in Ext\{A,B\}$ is called a quasi-selfadjoint extension of $\{A,B\}$ if
\be\lb{e31}
\dim(D(\tilde{A})/D(A))=\dim(D((\tilde{A})^*)/D(B)).
\ee

(2) A dual pair $\{A,B\}$ is called a correct dual pair if it admits a quasi-selfadjoint extension.
}
\end{defn}

Now, we shall give a necessary condition for quasi-selfadjoint extensions of dual pairs of linear  relations. Firstly, we give a property of self-adjoint extensions of linear relations.

The following result can be directly followed by Proposition 1.7.9 and Corollary 1.7.10 in \cite{Behrndt20}.

\bb{lem}\lb{l31}
Let $S\in LR(X)$ be Hermitian and $T\in LR(X)$ be a self-adjoint extension of $S$. Then
\be
n_{\pm}(S)=\dim(D(T)/D(S))=\dim(D(S^*)/D(T)).\nn
\ee
\end{lem}

\bb{thm}\lb{t1}
Let $A,B\in CR(X)$ and $\{A,B\}$ be a dual pair with $(A^*)_s|_{D(B)}=B_s$ and $(B^*)_s|_{D(A)}=A_s$. Suppose that $B^*(0)\cap N(A^*)=\{0\}$, $A^*(0)\cap N(B^*)=\{0\}$, and $\tilde{A}\in {\rm Ext}\{A,B\}$.
\begin{itemize}
\item[{\rm (i)}] The following formulae hold
\bea\lb{e32}
n(A^*,B)&=&\dim(D(B^*)/D(\tilde{A}))+\dim(D(A^*)/D((\tilde{A})^*))\nn\\
&=&\dim(D(\tilde{A})/D(A))+\dim(D((\tilde{A})^*)/D(B)).
\eea
\item[{\rm (ii)}] If $\tilde{A}$ is a  quasi-selfadjoint extension of $\{A,B\}$, then $n(A^*,B)$ is even and
\bea\lb{e33}
n(A^*,B)/2&=&\dim(D(B^*)/D(\tilde{A}))=\dim(D(A^*)/D((\tilde{A})^*))\nn\\
&=&\dim(D(\tilde{A})/D(A))=\dim(D((\tilde{A})^*)/D(B)).
\eea
\end{itemize}
\end{thm}

\Proof
(i) Let $T=\left(
             \begin{array}{cc}
               0 & \tilde{A} \\
               (\tilde{A})^* & 0 \\
             \end{array}
           \right)$ and $S$ be defined in \eqref{S}. By Lemma \ref{l4} and Proposition \ref{p2}, we have that $S$ is Hermitian and
\be\lb{e34}
S\subset T=\left(
             \begin{array}{cc}
               0 & \tilde{A} \\
               (\tilde{A})^* & 0 \\
             \end{array}
           \right)=T^*\subset S^*,
\ee
which implies that $T$ is a self-adjoint extension of $S$. Then
\be\lb{e35}
n_{\pm}(S)=\dim(D(T)/D(S))=\dim(D(S^*)/D(T))
\ee
by Lemma \ref{l31}. Note that
\be
D(T)=D((\tilde{A})^*)\times D(\tilde{A});D(S)=D(B)\times D(A);D(S^*)=D(A^*)\times D(B^*).\nn
\ee
Hence
\beaa
\dim (D(T)/D(S))&=&\dim(D(\tilde{A})/D(A))+\dim(D((\tilde{A})^*)/D(B));\\
\dim (D(S^*)/D(T))&=&\dim(D(B^*)/D(\tilde{A}))+\dim(D(A^*)/D((\tilde{A})^*))
\eeaa
which, together with \eqref{e12} and \eqref{e35}, implies that \eqref{e32} holds.

(ii) Suppose that $\tilde{A}$ is a  quasi-selfadjoint extension of $\{A,B\}$, i.e., \eqref{e31} holds. Note that
\bea\lb{e36}
n(A^*,B)&=&\dim(D(A^*)/D((\tilde{A})^*))+\dim(D((\tilde{A})^*)/D(B));\nn\\
n(B^*,A)&=&\dim(D(B^*)/D(\tilde{A}))+\dim(D(\tilde{A})/D(A)).
\eea
Then
\be
\dim(D(A^*)/D((\tilde{A})^*))=\dim(D(B^*)/D(\tilde{A})),\nn
\ee
by Proposition \ref{p1} and \eqref{e31}. This, together with \eqref{e31} and \eqref{e32}, implies that \eqref{e33} holds.
\qed

At the end of this section, we present a sufficient condition for a quasi-selfadjoint extension of a dual pair of linear relations.

\bb{thm}\lb{t2}
Let $A,B\in CR(X)$ and $\{A,B\}$ be a dual pair with $(A^*)_s|_{D(B)}=B_s$ and $(B^*)_s|_{D(A)}=A_s$. Suppose that $B^*(0)\cap N(A^*)=\{0\}$, $A^*(0)\cap N(B^*)=\{0\}$, and for some $\tilde{A}\in {\rm Ext}\{A,B\}$, we have
\be\lb{e37}
\dim(D(B^*)/D(\tilde{A}))=\dim(D(\tilde{A})/D(A)).
\ee
Then $\tilde{A}$ is a quasi-selfadjoint extension of $\{A,B\}$.
\end{thm}

\Proof
It follows from \eqref{e32} and \eqref{e37} that
\be
\dim(D(A^*)/D((\tilde{A})^*))=\dim(D((\tilde{A})^*))/D(B)),\nn
\ee
which together with \eqref{e36} and \eqref{e37} implies that
\beaa
n(A^*,B)&=&2\dim(D((\tilde{A})^*))/D(B));\\
n(B^*,A)&=&2\dim(D(\tilde{A})/D(A)).
\eeaa
Therefore,
\be
\dim(D(\tilde{A})/D(A))=\dim(D((\tilde{A})^*)/D(B)),\nn
\ee
i.e., $\tilde{A}$ is a quasi-selfadjoint extension of $\{A,B\}$.
\qed

\bb{rem}\lb{r2}
{\rm Theorems \ref{t1} and \ref{t2} are the generalizations of \cite[Propositions 2.10 and 2.11]{BrownK19} for dual pairs of operators to dual pairs of linear relations, respectively.}
\end{rem}

\section*{Statements and Declarations}
The authors have no competing interests to declare that are relevant to the content of this article.

\section*{Acknowledgements}

This work was supported by the National Natural Science Foundation of Shandong Province [Grants ZR2020MA012].

%\ifl
\section*{References}
%\fi

%%%
\end{document}